\documentclass[11pt]{article}
\usepackage{amsfonts}
\usepackage{color}

\usepackage{mathrsfs}
\usepackage[centertags]{amsmath}
\usepackage{amsfonts}
\usepackage{amssymb}
\usepackage{amsthm}
\usepackage{cases}
\usepackage{indentfirst}
\usepackage{epsfig}
\usepackage{graphicx}
\usepackage{color}
\usepackage{CJK}

\date{}

\definecolor{c20}{rgb}{0.,0.7,0.}
\definecolor{c30}{rgb}{0.,0.,1.}
\definecolor{c40}{rgb}{1,0.1,0.7}
\definecolor{c50}{rgb}{1,0,0}
\definecolor{c60}{rgb}{1,0.9,0.1}

\def\cL#1{\textcolor{c40}{#1}}
\def\cL#1{#1}
\def\ccL#1{\textcolor{c50}{#1}}
\def\ccL#1{#1}
\def\cE#1{#1}
\newcommand{\pk}[1]{\mathbb{P} \left( #1 \right)}

\newcommand{\expon}[1]{\exp\left(#1\right)}
\newcommand{\BQN}{\begin{eqnarray}}
\newcommand{\EQN}{\end{eqnarray}}
\newcommand{\BQNY}{\begin{eqnarray*}}
\newcommand{\EQNY}{\end{eqnarray*}}
\newcommand{\BS}{\begin{sat}}
\newcommand{\ES}{\end{sat}}
\newcommand{\BT}{\begin{theo}}
\newcommand{\ET}{\end{theo}}
\newcommand{\BL}{\begin{lem}}
\newcommand{\EL}{\end{lem}}
\newcommand{\BK}{\begin{korr}}
\newcommand{\EK}{\end{korr}}
\def\IF{\infty}

\newtheorem{theorem}{Theorem}[section]

\newtheorem{lemma}{Lemma}[section]

\newtheorem{remark}{Remark}[section]

\numberwithin{equation}{section}

\newcommand\abs[1]{\left\lvert #1 \right\rvert}
\def\P{\mathbb{P}}

\def\R{\operatorname*{\mathbb{R}}}
\def\N{\operatorname*{\mathbb{N}}}

\def\I#1{\mathbb{I}\{#1\}}
\allowdisplaybreaks[4]
\usepackage{setspace}
\begin{document}
\title{Higher-order expansions of  powered extremes of \\
normal samples }
\author{ {Wei Zhou\quad and\quad  Chengxiu Ling\footnote{Corresponding author. Email: lcx98@swu.edu.cn}}\\
{\small School of Mathematics and Statistics, Southwest University,
Chongqing, 400715, China}}
\date{}

\maketitle
\begin{quote}
{\bf Abstract.}~~ 
In this paper, higher-order expansions for distributions and densities of
powered extremes of standard normal random sequences are established under an optimal choice of normalized constants. Our findings refine the related results in Hall (1980). Furthermore, it is shown that the rate of convergence of  distributions/densities of normalized extremes depends in principle on the power index.
%Asymptotic expansions of the distributions of extremes for various random sequences have been received {much} attention  in both theoretical analysis and statistical applications.

{\bf Keywords.}~~Convergence rate; higher-order expansion; powered extremes.

{\bf MSC:}~~ Primary 62E20; 60E05; secondary 60F15; 60G15
\end{quote}

\section{Introduction}
\label{sec1}
Let $X_1, X_2, \ldots$ be independent random variables with common \cL{standard normal distribution function (df)} $\Phi(x)$, and let $M_n=\max(X_1,X_2,\ldots,X_n)$.
It is well-known that $\Phi(x)$ belongs to the max-domain of attraction of Gumbel extreme value distribution $\Lambda(x)= \expon{-e^{-x}}$,  i.e.,
there exist norming
constants $a_{n}>0$ and $b_{n}\in\R$ such that (see e.g., \ccL{Embrechts et al. (1997)})
\begin{eqnarray}
\label{eq1.1}
\lim_{n\to\IF}\P\left(M_n \leq {a_n}x+b_n\right) = \Lambda(x), \quad x\in\R.
\end{eqnarray}
%Hall (1979) showed that $1/\log n$ is the best convergence rate of \eqref{eq1.1}. %For more related work, see  and Nair (1981). %It was shown that the convergence rate of \eqref{eq1.1}
%is proportional to $1/\log n$ with optimal norming constants $a_n$ and $b_n$.
Further, %\cL{motivated by Haldance and Jayakar (1963),} 
Hall (1980) established the
asymptotic distribution behavior of  normalized
$\abs{M_n}^t$, the powered extremes for given power index $t>0$. Precisely speaking, with %$\phi(x)$ the probability density function of an $N(0,1)$ random variable, and
$b_{n}>0$ the solution of the following equation (see e.g., Gasulla et al. (2015) for related discussions on $b_n$)  
\begin{eqnarray}
\label{eq1.2}
2\pi b_n^2\exp\left(b_n^2\right)=n^2, \quad \forall n\in\N,
\end{eqnarray}
Hall (1980) showed that
\begin{eqnarray}
\label{eq1.4}
\lim_{n \to \infty} b_n^{2+2\I{t=2}} \Big(\P\left(\abs{M_n}^t\leq
c_n x+{d_n}\right)-\Lambda(x)\Big)
=\Lambda(x)e^{-x} k_1(t, x),\quad x\in\R,
\end{eqnarray}
where, with $\I{\cdot}$ the indicator function,
\begin{eqnarray}
\label{eq1.5}
&&
c_n=t b_n^{t-2} - 2b_n^{-2}\I{t=2}, \quad d_n=b_n^t - 2b_n^{-2}\I{t=2},\\
&&k_1(t, x)=\left(1+x+\frac{2-t}2 x^2\right)\I{t\neq2} - \left(\frac72+3x+x^2\right)\I{t=2}
\label{eq1.3}.
\end{eqnarray}
 Recently, Liao et al. (2013a) and Hashorva
et al. (2014) studied respectively expansions of asymptotic distributions of normalized extremes for \ccL{logarithmic}  skew-normal
distributions and bivariate normal triangular arrays. \cL{Liu and Liu (2013) and Chen and Huang (2014) studied respectively  uniform convergence rates of distributions of normalized extremes for Maxwell samples and asymmetry normal samples. 
Other related work on extreme value distributions, densities and \ccL{moments} of given distributions and
their associated uniform convergence rates, we refer to %Nair (1981), Smith (1982), Omey (1988), 
de Haan and Resnick (1996), Withers and Nadarajah (2010), Liao et al. (2013b), Liao et al. (2014), Peng et al. (2014) and \ccL{Li and Li (2015)}, and references therein.} 

In this paper, we aim to establish \ccL{higher-order} expansions
of distributions and densities of powered extremes
for normal random samples. The motivation is two-folded. One comes
from the importance of determining the efficiency of distribution/density approximations to its Gumbel limit law as shown by the contributions mentioned above. The other is that the powered normal
laws are one challenging alternative to normal ones due to its mathematical properties including skewness, heavy tails etc.

\cL{The contribution of this paper is to establish the rate of convergence of \eqref{eq1.4}. Moreover, we find out that} higher-order expansions of densities of powered extremes \ccL{display} similar asymptotic structures as those for \ccL{their} distributions.

The rest of the paper is organized as follows. Section \ref{sec2}
presents main results. All the proofs are relegated to Section \ref{sec3}.

\section{Main Results}
\label{sec2}
In this section, we shall establish higher-order expansions of distributions and densities of powered extremes under   normalization (see Theorems \ref{the2.0} and \ref{the2.1} below). It shows that the convergence rates are different \cL{between the two cases that} the \ccL{power} index $t=2$ and $0<t\neq2$. 

\cL{In the sequel, we shall} keep the notation given in Section \ref{sec1} unless stated otherwise.
Further, let \cL{for $x\in\R$ and $t>0$}
\begin{eqnarray} \label{def_k2}
k_2(t,x)=\left\{
\begin{array}{ll}
\frac{43}{3}+14x+6x^2+\frac{4}{3}x^3, & t=2;\\
- \left(3+3x+\frac{3}{2}x^2 + {\frac{(2-t)(2t+1)}{6}x^3}+\frac{(t-2)^2}{8}x^4  - \frac{e^{-x}}{2}\left(1+x+\frac{2-t}{2}x^2\right)^2 \right), & t\neq2
\end{array}
\right.
\end{eqnarray}
and
\BQN \label{def_kappa}
&&\varpi(t, x)= \left\{
\begin{array}{ll}
 \displaystyle \frac{1}{2}+x+x^2 - e^{-x}\left(\frac{7}{2}+3x+x^2\right), & t=2; \\
\displaystyle x \left(1-t + \frac{t-2}2 x\right) + e^{-x}\left(1+x+\frac{2-t}{2}x^2\right) , & t\neq2
\end{array}
\right. \\
\label{def_tau}
&& \tau(t, x) =\I{t=2}\left( e^{-x}\left(\frac{43}{3}+14x+6x^2+\frac{4}{3}x^3\right)  - \left(\frac{1}{3}+2x+2x^2+\frac{4}{3}x^3\right)\right) \nonumber \\
&&\qquad \qquad+\I{t\neq2} \left( x e^{-x}
\left(1-t + \frac{t-2}2 x\right) \left(1+x+\frac{2-t}{2}x^2\right) \right. \nonumber \\
&&\qquad \qquad+x^2\left( \frac{(1-t)(1-2t)}{2} + {\frac{5(1-t)(t-2)}{6} }x + \frac{(t-2)^2}8 x^2  \right)- e^{-x} \Bigg(3+3x \nonumber  \\
&&\qquad \qquad\left.\left.+\frac{3}{2}x^2 + {\frac{(2-t)(2t+1)}{6}x^3}+\frac{(t-2)^2}{8}x^4  - \frac{e^{-x}}{2}\left(1+x+\frac{2-t}{2}x^2\right)^2 \right) \right). \quad
\EQN
\begin{theorem}
\label{the2.0} \cL{Let $M_n=\max(X_1, \ldots, X_n)$ with $\{X_n, n\ge1\}$ a sequence of independent random variables with common df $\Phi(x)$.  Let further} $b_n, c_n$ and $d_n$ are those given by \eqref{eq1.2} and \eqref{eq1.5}, respectively. We have for any $t>0$ and as $n\to\IF$
\begin{eqnarray*}
\pk{\abs{M_n}^t\le c_n x+ d_n}= \Lambda(x) +
 b_n^{-2-2\I{t=2}}\Lambda'(x)\left(k_1(t, x) + b_n^{-2} k_2(t, x) +
O(b_n^{-4}) \right), \quad x\in\R,
\end{eqnarray*}
where $\Lambda'(x)= \Lambda(x) e^{-x}$, and $k_1(t,
x)$ and $k_2(t,x)$ are {those}
 given by \eqref{eq1.3} and
\eqref{def_k2}, respectively.
\end{theorem}
\begin{remark}
Theorem \ref{the2.0} gives the accurate convergence rate of
\eqref{eq1.4},  which is proportional to $1/\log n$ for all $t>0$
since $b^{2}_{n}{\sim} 2\log n$ as $n\to\infty$ due to \eqref{eq1.2}.
\end{remark}
Next, we shall establish the rate of convergence of density approximation of 
$(|M_n|^t-d_n)/c_n$ to Gumbel extreme value density function. 
\begin{theorem}
\label{the2.1} Under the assumptions of Theorem \ref{the2.0}, we
have as $n\to\IF$
\begin{eqnarray*}
\frac d{dx}\pk{\abs{M_n}^t\le c_n x+ d_n} = \Lambda'(x)\left( 1+ b_n^{-2-2\I{t=2}}\left(
\varpi(t, x) + b_n^{-2} \tau(t, x) + O(b_n^{-4}) \right) \right),
\ \  x\in\R,
\end{eqnarray*}
where $\varpi(t, x)$ and $\tau(t, x)$ are {those} given by
\eqref{def_kappa} and \eqref{def_tau}, respectively.
\end{theorem}
\begin{remark}
(1) %We see that the density of $(|M_{n}|^{t}-d_{n})/c_{n}$ converges to that of Gumbel extreme value distribution $\Lambda(x)$, and 
We see that higher-order approximations of densities of $(|M_{n}|^{t}-d_{n})/c_{n}$ \cL{to Gumbel extreme value density function} possess the same structure as those for their distributions, i.e., the  second-order
convergence rate is rather faster for $t=2$ than that for the other cases,  while the third-order one is the same $1/\log n$ for all $t>0$. \\
%\cL{(2) One might apply Theorem \ref{the2.1} to derive higher-order expansions of moments of $(|M_n|^t-d_n)/c_n$. \ccL{For related work}, see Theorem 3.1 and Theorem 2.3 respectively by Nair (1981) and Li and Li (2015) for normal distributions and general error distributions.} \\
(2) It might be possible to investigate the expansions under consideration for powered $k$th extremes of
normal distributions following the similar arguments, see Theorem 1 by Hall (1980).
\end{remark}

\section{Proofs}
\label{sec3} In this section, we shall present the proofs of Theorems \ref{the2.0} and \ref{the2.1}. To this end, we first establish three lemmas. The first one is concerned with a crucial
expansion of the {survival function $1-\Phi(x)$}. Lemmas \ref{lemma2} and
\ref{lemma3} specify the expansions of the two terms of densities of
$(\abs{M_n}^t-d_n)/c_n$ (see \eqref{exp_density} below). 
Hereafter, all the limit relations are for $n\to\IF$ unless otherwise stated. 
\begin{lemma}
\label{lemma1} Let $\phi(x)$ be the probability density function of an $N(0,1)$ random variable. We have for any given $L\in\N$
\begin{equation*}
 1-\Phi(x) =
x^{-1}\phi(x) \left( \sum_{k=0}^L (-1)^k(2k+1)!!x^{-2k} +
O(x^{-2(L+1)}) \right)
\end{equation*}
as $ x\to\IF$, where $(2k+1)!!=(2k)!/\cE{(2^{k}k!)}$ with
$m!=m(m-1)(m-2)\cdots 2\cdot 1$.
\end{lemma}
\noindent \textbf{Proof.} Note that $1-\Phi(x) = x^{-1}\phi(x) \int_0^\IF e^{-t}e^{-(t/x)^2/2}\ dt,\; x>0$. Therefore, the claim follows by applying the dominated convergence theorem  to integrate term by term the \ccL{$L$th} Taylor's expansion of $f(u)=e^{-u},\; u=(t/x)^2/2$ at $u=0$ (recall $x\to\IF$).  \qed
\begin{lemma}
\label{lemma2} Let  \cL{$\ccL{\nu_n(x, t)}= \Phi^{n-1}\left((c_nx+{d_n})^{1/t}\right)-
\left(1-\Phi\left((c_nx+{d_n})^{1/t}\right)\right)^{n-1},\; x\in\R, \; t>0$} with $c_{n}$ and $d_{n}$ given by \eqref{eq1.5}. \\
(1) For $0<t\neq2$, we have as $n\to\IF$
\begin{eqnarray}
\label{eq3.2}
& & \nu_n(x, t) = \Lambda(x) \left(1+\frac{e^{-x}}{b_n^{2}}\left(1+x+\frac{2-t}{2}x^2\right) \right.- \frac{e^{-x}}{b_n^4}\left(3+3x+\frac{3}{2}x^2\right.\nonumber
\\
&&\quad \left.\left. + {\frac{(2-t)(2t+1)}{6}x^3}+\frac{(t-2)^2}{8}x^4  - \frac{e^{-x}}{2}\left(1+x+\frac{2-t}{2}x^2\right)^2 \right)+O(b_n^{-6})\right).
\end{eqnarray}
(2) For $t=2$, we have as $n\to\IF$
\begin{eqnarray}
\label{eq3.3}
\nu_n(x, t)= \Lambda(x)\left(1- \frac{e^{-x}}{b_n^4}\left(\frac{7}{2}+3x+x^2\right) +  \frac{e^{-x}}{b_n^6}\left(\frac{43}{3}+14x+6x^2+\frac{4}{3}x^3\right)+O(b_n^{-8})\right).
\end{eqnarray}
\end{lemma}

\noindent \textbf{Proof.}~~\cL{Note for fixed $x\in\R$ and large $n$ that $c_n x+d_n>0$. % with $c_n$ and $d_n$ given  by \eqref{eq1.5}. 
Set below} $g_n(t)=(c_n x + d_n)^{1/t}$.

\underline{(1) For $0<t\neq 2$.} Using the following
Taylor's expansion
$$(1+x)^\alpha = 1+\alpha x + \frac{\alpha(\alpha-1)} 2 x^2 + \frac{\alpha(\alpha-1)(\alpha-2)} 6x^3 (1+O(x)), \quad x\to
0,\ \alpha\in\R$$  
\cL{and the fact by \eqref{eq1.2} that  $b_n^{2}\sim \ccL{2}\log n$ for large $n$,} we have by \eqref{eq1.5}
\BQN \label{eq3.1} g_n^\alpha(t)
=b_n^\alpha \left( 1+ \frac{\alpha x}{b_n^2} + \frac{\alpha
(\alpha - t)}{2b_n^4} x^2 + \frac{\alpha (\alpha - t) (\alpha
-2t)}{6b_n^6} x^3 (1+ O(b_n^{-2}))\right). 
\EQN
Applying \eqref{eq3.1} with $\alpha=-1$ and $\alpha=2$, we have
\BQN\label{eq3.phi} \ccL{\frac{\phi(g_n(t))}{g_n(t)}} 
 &=&
b_n^{-1}\left( 1 - \frac x {b_n^2} + \frac{1+t}{2b_n^4} x^2+ O(b_n^{-6}) \right)
\nonumber
\\
&&
\times  \frac 1{\sqrt{2\pi}} \expon{-\frac{b_n^2}2 \left(1+\frac 2 {b_n^2} x + \frac{2-t}{b_n^4} x^2 +\frac{2(t-1)(t-2)}{3b_n^6}x^3 + O(b_n^{-8})\right)}
\\
&=&
\frac{\phi(b_n)}{b_n} e^{-x} \left( 1 - \frac x {b_n^2} + \frac{1+t}{2b_n^4} x^2 +  O(b_n^{-6}) \right)\nonumber \\
&& \times 
\left(1 + \frac{t-2}{2b_n^2}x^2 - \frac{(t-1)(t-2)}{3b_n^4}x^3  + \frac{(t-2)^2}{8b_n^4}x^4 +  O(b_n^{-6}) \right)
\nonumber
\\
&=& \frac{e^{-x}}n \left(1 +  \ccL{\frac x{b_n^2}} \left(\frac{t-2}2 x -1\right) \right. \nonumber \\ 
&& \left. +\frac{x^2}{b_n^4}\left(\frac{t+1}2 -\frac{(t-2)(2t+1)}6 x + \frac{(t-2)^2}8x^2 \right) +  O(b_n^{-6})\right),
\nonumber
\EQN
where the second equality \cL{holds} since $e^x=1+x+x^2/2(1+O(x)), \ x\to0$ and the third equality \cL{follows by \eqref{eq1.2}.}
Further, it follows by \eqref{eq3.1} with $\alpha=-2$ and $\alpha=-4$ that
\BQNY
1-g_n^{-2}(t) + 3g_n^{-4}(t) + O(g_n^{-6}(t))
&=& 1- b_n^{-2}\left(1-2 b_n^{-2}x+ O(b_n^{-2})\right) + 3b_n^{-4}(1+O(b_n^{-2}))
\\
&=& 1- b_n^{-2} +b_n^{-4}(2x+3)+ O(b_n^{-6}).
\EQNY
Therefore, it follows further by Lemma \ref{lemma1} with $L=2$ that
\begin{eqnarray}
\label{def:coef}
 1-\Phi\left((g_n(t)\right)
& \cL{=} & \frac{\phi(g_n(t))}{g_n(t)} \big( 1-g_n^{-2}(t) + 3g_n^{-4}(t) + O(g_n^{-6}(t)) \big) \nonumber
\\
& = &
\frac{e^{-x}}n \left(1- \frac1{b_n^{2}}\left(1+x+\frac{2-t}{2}x^2\right) \right. \nonumber \\
&& +\left.
 \frac1{b_n^{4}}\left(3+3x+\frac{3}{2}x^2+\frac{(2-t)(2t+1)}{6}x^3+\frac{(t-2)^2}{8}x^4\right)+O(b_n^{-6})\right)
 \nonumber\\
 &=:& n^{-1}e^{-x} \left(1-\vartheta_1b_n^{-2}+\vartheta_2b_n^{-4}+O(b_n^{-6})\right),
\end{eqnarray}
which together with the fact that $b_n^2 \sim 2\log n$ and
$\log(1-x)=-x(1+{O}(x)),\ x\to0$ implies that
\begin{eqnarray}
\label{eq3.4}
\Phi^{n-1}\left((g_n(t) \right) &=&  \exp\big((n-1)\log(1-(1-\Phi \left(g_n(t) \right))) \big)
\nonumber
\\
&=&
\expon{ -e^{-x} \left(1 - \vartheta_1b_n^{-2}+\vartheta_2b_n^{-4}+O(b_n^{-6})\right)}
\nonumber  \\
&=&
\Lambda(x)\left( 1+\frac{\vartheta_1e^{-x}}{b_n^2} +\frac{e^{-x}}{b_n^4}
\left(\frac{\vartheta_1^2}2e^{-x}-
\vartheta_2\right)
 +O(b_n^{-6}) \right)
\end{eqnarray}
and
\BQN \label{eq.Ompart}
\left(1-\Phi\left(g_n(t)\right)\right)^{n-1} = (n^{-1}e^{-x}(1+O(b_n^{-2})))^{n-1} = o(b_n^{-\alpha}), \quad \alpha \ge 6.
\EQN
Hence \eqref{eq3.2} follows by recalling $\vartheta_1$ and $\vartheta_2$ given by \eqref{def:coef}.

\underline{(2) For $t=2$.} We shall verify \eqref{eq3.3} by similar
arguments to {those for} the case $0<t\neq2$.  Recalling $c_n =
2(1-\cE{b}_n^{-2}), d_n = b_n^2-2b_n^{-2}$ for $t=2$, we have $$h_n:=
g_n(2)= (c_nx+d_n)^{1/2}=b_n(1+2xb_n^{-2}-2(1+x)b_n^{-4})^{1/2}.$$
Hence,
\begin{eqnarray*}
h_n^{\alpha} = b_n^\alpha \left(1+\frac{\alpha x}{b_n^2} -\frac\alpha{b_n^4}\left(1+x+\frac{2-\alpha}2x^2\right)
+\frac{\alpha(2-\alpha)x}{b_n^6}\left(1+x+\frac{4-\alpha}6x^2\right) +  O(b_n^{-8})\right).
\end{eqnarray*}
Using $e^x=1+x+x^2/2+x^3/6 + O(x^4),\ x\to0$ and the above equality with $\alpha=-1$, we have
\begin{eqnarray} \label{eq0}
\frac{\phi(h_n)}{h_n}
&=&
\frac{\phi(b_n)}{b_n} e^{-x}\left( 1+\frac{1+x}{b_n^2}+\frac{(1+x)^2}{2b_n^4}+\frac{(1+x)^3}{6b_n^{6}}+ O(b_n^{-8})\right) \nonumber\\
&&\times \left( 1 - \frac x{b_n^2} + \frac1{b_n^4}\left(1+x+\frac32 x^2\right) - \frac{3x}{b_n^6}\left(1+x+\frac56 x^2\right)+ O(b_n^{-8})\right) \nonumber
\\
&=& \frac{e^{-x}}n \left(1 + \frac1{b_n^2}+\frac1{b_n^4} \left(x^2+x+\frac32\right) - \frac1{b_n^6}\left( \frac43x^3 + x^2 + x - \frac76\right) + O(b_n^{-8})\right).\qquad
\end{eqnarray}
Additionally,
\begin{eqnarray*}
& & 1-h_n^{-2}+3h_n^{-4}-15h_n^{-6} + O(b_n^{-8})
\nonumber
\\
&&=
1- b_n^{-2}\Big(1-2xb_n^{-2}+2(1+x+2x^2)b_n^{-4}-8x(1+x+x^2)b_n^{-6}+O(b_n^{-8})\Big)
\nonumber
\\
& &\quad +3b_n^{-4}\Big(1-4xb_n^{-2}+O(b_n^{-2})\Big)
 -15b_n^{-6}\Big(1+O(b_n^{-2})\Big) + O(b_n^{-8})
\\
&&=
1-b_n^{-2} + b_n^{-4}(2x+3) - b_n^{-6}(17+14x+4x^2) + O(b_n^{-8}).\quad
\end{eqnarray*}
Consequently, a straightforward application of Lemma \ref{lemma1} with $L=3$ yields that
\begin{eqnarray*}
1- \Phi(h_n) = \frac{e^{-x}}n\left(1 + \frac1{b_n^4}\left(\frac{7}{2}+3x+x^2\right)- \frac1{b_n^6}\left(\frac{43}{3}+14x+6x^2+\frac{4}{3}x^3\right) + O(b_n^{-8})\right).
\end{eqnarray*}
The rest proof follows by the same arguments as for \eqref{eq3.4} and \eqref{eq.Ompart} with $g_n(t)$ replaced by $h_n$.  We complete the proof of Lemma \ref{lemma2}.
 \qed

\begin{lemma}\label{lemma3} Let $c_{n}$ and $d_{n}$ be given by \eqref{eq1.5}. We have as $n\to\IF$
\begin{eqnarray*}
&&n\frac d{dx}\Phi((c_nx + d_n)^{1/t})
\nonumber \\
&&= \left\{ \begin{array}{ll}
e^{-x} \left( 1+\frac x{b_n^2} \left(1-t + \frac{t-2}2 x\right) +\frac{x^2}{b_n^4} \left( \frac{(1-t)(1-2t)}{2} + {\frac{5(1-t)(t-2)}{6} }x + \frac{(t-2)^2}8 x^2  \right) + O(b_n^{-6}) \right), &t\neq2;\\
 e^{-x} \left(1+ \ccL{b_n^{-4}}\left( \frac{1}{2}+x+x^2\right) - b_n^{-6}\left(\frac{1}{3}+2x+2x^2+\frac{4}{3}x^3\right)+O(b_n^{-8})\right), & t=2.
\end{array}
\right. 
\end{eqnarray*}
\end{lemma}

\noindent \textbf{Proof.}~~Clearly, we have 
\BQNY
n\frac d{dx}\Phi((c_nx + d_n)^{1/t}) =t^{-1}nc_n\left(c_nx+{d_n}\right)^{\ccL{1/t}-1}\phi\left((c_nx+{d_n})^{1/t}\right).
\EQNY 
For $0<t\neq2$, we have by \eqref{eq1.5} that $c_n=tb_n^{t-2}$.  It follows by \eqref{eq3.1} with $\alpha=1-t$, and \eqref{eq3.phi} for the expansion of $\phi(g_n(t))$ with $g_n(t)=(c_nx+{d_n})^{1/t}$ that
\begin{eqnarray*}
n\frac d{dx}\Phi(g_n(t))&=& \frac n{b_n} \left( 1+ \frac{(1-t)x }{b_n^2}+\frac{(1-t)(1-2t)}{2b_n^4}x^2+O(b_n^{-6}) \right) \\
&& \times \frac 1{\sqrt{2\pi}} \expon{-\frac12b_n^2 \left(1+\frac 2 {b_n^2} x + \frac{2-t}{b_n^4} x^2 +\frac{2(t-1)(t-2)}{3b_n^6}x^3 + O(b_n^{-8})\right)}
\\
&=& e^{-x} \left( 1+ \frac{(1-t)x }{b_n^2}+\frac{(1-t)(1-2t)}{2b_n^4}x^2+O(b_n^{-6}) \right) \\
&& \times \left(1 + \frac{t-2}{2b_n^2}x^2 - \frac{(t-1)(t-2)}{3b_n^4}x^3  + \frac{(t-2)^2}{8b_n^4}x^4 +  O(b_n^{-6}) \right) \\
&=& e^{-x} \left( 1+\frac x{b_n^2} \left(1-t + \frac{t-2}2 x\right) \right.\\
&& + \left.\frac{x^2}{b_n^4} \left( \frac{(1-t)(1-2t)}{2} + \frac{5(1-t)(t-2)}{6} x + \frac{(t-2)^2}8 x^2\right) +  O(b_n^{-6})
\right). 
\end{eqnarray*}
%We completes the proof of \eqref{eq3.9} for $0<t\neq2$.
\cL{For $t=2$,} recalling that $c_n= 2(1-b_n^{-2})$, it follows by \eqref{eq0} that,
%the left-hand side of \eqref{eq3.9} equals
\begin{eqnarray*}
n\frac d{dx}\Phi(g_n(t))&=& e^{-x} \left(1-\frac1{b_n^2}\right) \left(1 + \frac1{b_n^2}+ \frac1{b_n^4} \left(x^2+x+\frac32\right) - \frac1{b_n^6}\left( \frac43x^3 + x^2 + x - \frac76\right) + O(b_n^{-8})\right)
\\
&=& e^{-x} \left(1+ \frac1{b_n^4}\left( \frac{1}{2}+x+x^2\right) - \frac1{b_n^6}\left(\frac{1}{3}+2x+2x^2+\frac{4}{3}x^3\right)+O(b_n^{-8})\right).
\end{eqnarray*}
The proof is complete.  \qed

\noindent \textbf{Proof of Theorem 2.1}~~Note that {for fixed $x\in\R$ and large $n$, $c_n x+d_n>0$ and} the distribution function of $(\abs{M_n}^t - d_n)/c_n$ is as follows.
\BQN
\label{eq_dis}
\pk{\abs{M_n}^t \le c_n x + d_n} = \Phi^n\left((c_nx+{d_n})^{1/t}\right)-
\left(1-\Phi\left((c_nx+{d_n})^{1/t}\right)\right)^n.
\EQN
For $0<t\neq2$, it follows by \eqref{def:coef} and similar arguments as for \eqref{eq3.4} and \eqref{eq.Ompart} that
\begin{eqnarray*}
\Phi^{n}\left((g_n(t) \right) =
\Lambda(x)\left( 1+\frac{\vartheta_1e^{-x}}{b_n^2} +\frac{e^{-x}}{b_n^4}
\left(\frac{\vartheta_1^2}2e^{-x}-
\vartheta_2\right)
 +O(b_n^{-6}) \right),
\end{eqnarray*}
where $\vartheta_1$ and $\vartheta_2$ are given by \eqref{def:coef}, and
\BQNY
\left(1-\Phi\left(g_n(t)\right)\right)^{n} = o(b_n^{-\alpha}), \quad \alpha \ge 6.
\EQNY
 Therefore, \eqref{eq3.2} holds \cL{with $\nu_n(x, t)$ replaced by the distribution of $(\abs{M_n}^t - d_n)/c_n$ (see \eqref{eq_dis} above).}  
 
\cL{Similarly, for the case $t=2$, we have \eqref{eq3.3} holds with \ccL{$\nu_n(x, t)$} replaced by the distribution of $(\abs{M_n}^t - d_n)/c_n$.  Consequently, the desired result in Theorem \ref{the2.0} is obtained for all $t>0$.} \qed

\noindent \textbf{Proof of Theorem 2.2.}~~By the symmetry of standard normal distribution $\Phi(-x)=1-\Phi(x)$, we have
\BQN
\label{exp_density}
\lefteqn{\frac d{dx}\pk{\abs{M_n}^t\le c_n x+ d_n} = \ccL{n\frac d{dx}\Phi((c_nx + d_n)^{1/t})}} \nonumber \\
&&\times
\Big(\Phi^{n-1}\left((c_nx+{d_n})^{1/t}\right)-
\left(1-\Phi\left((c_nx+{d_n})^{1/t}\right)\right)^{n-1}\Big).
\EQN
Hence, for $0<t \neq2$, %applying  \eqref{eq3.2} and \eqref{eq3.9} for $0<t \neq2$, we have
it follows by Lemmas \ref{lemma2} and \ref{lemma3} that 
\begin{eqnarray*}
& & \frac{1}{\Lambda'(x)}\frac d{dx}\pk{\abs{M_n}^t\le c_n x+ d_n} -1
= \left( 1+\frac x{b_n^2} \left(1-t + \frac{t-2}2 x\right) \right. \\
&& \quad \left.+\frac{x^2}{b_n^4} \left( \frac{(1-t)(1-2t)}{2} + {\frac{5(1-t)(t-2)}{6} }x + \frac{(t-2)^2}8 x^2  \right) + O(b_n^{-6})
\right)
\\
&& \quad \times \left(1+\frac{e^{-x}}{b_n^{2}}\left(1+x+\frac{2-t}{2}x^2\right) - \frac{e^{-x}}{b_n^4}\right.\left(3+3x+\frac{3}{2}x^2 \right.
\\
&& \quad\left. +\left. {\frac{(2-t)(2t+1)}{6}x^3}+\frac{(t-2)^2}{8}x^4  - \frac{e^{-x}}{2}\left(1+x+\frac{2-t}{2}x^2\right)^2 \right)+O(b_n^{-6})\right) - 1\\
&&= \frac{1} {b_n^2} \left(x \left(1-t + \frac{t-2}2 x\right) + e^{-x}\left(1+x+\frac{2-t}{2}x^2\right) \right) \\
&&\quad +\frac{1} {b_n^4}\bigg( x e^{-x}
\left(1-t + \frac{t-2}2 x\right) \left(1+x+\frac{2-t}{2}x^2\right) \\
&&\quad+x^2\left( \frac{(1-t)(1-2t)}{2} + {\frac{5(1-t)(t-2)}{6} }x + \frac{(t-2)^2}8 x^2  \right)- e^{-x}  \\
&&\quad \left.\times\left(3+3x+\frac{3}{2}x^2 + {\frac{(2-t)(2t+1)}{6}x^3}+\frac{(t-2)^2}{8}x^4  - \frac{e^{-x}}{2}\left(1+x+\frac{2-t}{2}x^2\right)^2 \right) \right) + O(b_n^{-6}) \\
&&= b_n^{-2} \varpi(t,x)+ b_n^{-4}\tau(t,x)+
O(b_n^{-6}),
\end{eqnarray*}
\ccL{which completes the proof of} Theorem \ref{the2.1} for $t\neq2$.  Here $\varpi(t,x)$
and $\tau(t,x)$ are given by \eqref{def_kappa} and \eqref{def_tau},
respectively.

Next, for $t=2$, using  Lemmas \ref{lemma2} and \ref{lemma3} \ccL{together with} \eqref{exp_density}, we have
\begin{eqnarray*}
& & \frac{1}{\Lambda'(x)}\frac d{dx}\pk{\abs{M_n}^t\le c_n x+ d_n}-1 \\
&& = \left(1- \frac{e^{-x}}{b_n^4}\left(\frac{7}{2}+3x+x^2\right) +  \frac{e^{-x}}{b_n^6}\left(\frac{43}{3}+14x+6x^2+\frac{4}{3}x^3\right)+O(b_n^{-8})\right)\\
& &\quad \times \left(1+\frac1{b_n^4}\left( \frac{1}{2}+x+x^2\right) - \frac1{b_n^6}\left(\frac{1}{3}+2x+2x^2+\frac{4}{3}x^3\right)+O(b_n^{-8})\right) - 1\\
&& =\frac1{b_n^4} \left(  \frac{1}{2}+x+x^2 - e^{-x}\left(\frac{7}{2}+3x+x^2\right) \right) \\
&&\quad + \frac1{b_n^6}
\left(e^{-x}\left(\frac{43}{3}+14x+6x^2+\frac{4}{3}x^3\right)  - \left(\frac{1}{3}+2x+2x^2+\frac{4}{3}x^3\right) \right) + O(b_n^{-8})
\end{eqnarray*}
deriving the desired results. We complete the proof of Theorem \ref{the2.1}. \qed

\vspace{1cm}

\noindent{\bf Acknowledgements}~~The authors would like to thank \ccL{the Editor-in-Chief} and the referees for careful reading and comments which greatly improved the paper.
The work was supported by the National Natural Science Foundation of China (grant No.11171275) and Fundamental Research Funds for the Central Universities (XDJK2016C118, SWU115089).

\end{document}